\documentclass[12pt]{article}

\usepackage{amsthm,amsmath,amssymb,latexsym,url}

\newtheorem{thm}{Theorem}[section]

\newtheorem{cor}[thm]{Corollary}

\newtheorem{lem}[thm]{Lemma}
\newtheorem{conj}[thm]{Conjecture}
\newtheorem{problem}[thm]{Problem}

\theoremstyle{remark}
\newtheorem*{rmk}{Remark}
\def\Z{\mathbb{Z}}
\def\N{\mathbb{N}}

\def\pf{\noindent {\it Proof.} }

\numberwithin{equation}{section}
\renewcommand{\qed}{{\hfill\rule{4pt}{7pt}}\medskip}

\begin{document}
\baselineskip=17pt
\begin{center}
{\Large\bf  Factors of alternating sums of products of
binomial and $q$-binomial coefficients}
\end{center}
\vskip 2mm
\centerline{Victor J. W. Guo, Fr\'ed\'eric Jouhet  and Jiang Zeng}

\vskip 0.7cm {\small \noindent{\bf Abstract.}
In this paper we study the factors of some alternating sums of products of
binomial and $q$-binomial coefficients. We prove that
for all positive integers
$n_1, \ldots, n_m$, $n_{m+1}=n_1$, and $0\leq j\leq m-1$,
$${n_1+n_{m}\brack n_1}^{-1}\sum_{k=-n_1}^{n_1}(-1)^kq^{jk^2+{k\choose 2}}
\prod_{i=1}^m {n_i+n_{i+1}\brack n_i+k}\in \N[q],$$ which
generalizes a result of Calkin [Acta Arith. 86 (1998), 17--26].
Moreover, we show that for all positive integers $n$, $r$ and $j$,
$$
{2n\brack n}^{-1}{2j\brack j}
\sum_{k=j}^n(-1)^{n-k}q^{A}\frac{1-q^{2k+1}}{1-q^{n+k+1}}
{2n\brack n-k}{k+j\brack k-j}^r\in \N[q],
$$
where $A=(r-1){n\choose 2}+r{j+1\choose 2}+{k\choose 2}-rjk$,
which solves a problem raised by Zudilin
[Electron. J. Combin. 11 (2004), \#R22].


\section{Introduction}
In 1998, Calkin~\cite{Calkin} proved that for all positive integers $m$ and $n$,
\begin{align}
{2n\choose n}^{-1}\sum_{k=-n}^n (-1)^k{2n\choose n+k}^m \label{eq:Calkin}
\end{align}
is an integer by arithmetical techniques. For $m=1,2$ and 3,
by the binomial theorem, Kummer's formula and Dixon's formula,
it is easy to see that \eqref{eq:Calkin} is equal to
$0$, $1$ and ${3n\choose n}$, respectively.
Recently in the study of finite forms of the
Rogers-Ramanujan identities~\cite{GJZ} we stumbled across
\eqref{eq:Calkin} for $m=4$ and $m=5$, which gives
\begin{align*}
\sum_{k=0}^n {2n+k\choose k}{2n\choose n+k}^2
\quad\text{and}\quad
\sum_{k=0}^n {3n-k\choose n-k}{2n+k\choose k}{2n\choose n+k}^2,
\end{align*}
respectively. Indeed, de Bruijn~\cite{Bruijn} has shown that
for $m\geq 4$ there is  no closed form for \eqref{eq:Calkin}
by asymptotic techniques.
Our first objective is to give a $q$-analogue of Calkin's
result, which also implies that \eqref{eq:Calkin} is positive for $m\geq 2$.

In 2004, Zudilin~\cite{Zu} proved that for all positive integers $n$, $j$ and $r$,
\begin{align}\label{eq:strehl}
{2n\choose n}^{-1}{2j\choose j}\sum_{k=j}^n (-1)^{n-k}
\frac{2k+1}{n+k+1}{2n\choose n-k}{k+j\choose k-j}^r\in \Z,
\end{align}
which was originally observed
by Strehl~\cite{St} in 1994. In fact, Zudilin's motivation was to solve
the following problem,  which was   raised by
Schmidt~\cite{Sc} in 1992 and was apparently not related to Calkin's result.
\begin{problem}[Schmidt~\cite{Sc}]\label{prob:Sch}
For any integer $r\geq 2$, define a sequence of numbers $\{c_k^{(r)}\}_{k\in \N}$,
independent of the parameter $n$, by
\begin{equation*}
\sum_{k=0}^n{n\choose k}^r{n+k\choose k}^r
=\sum_{k=0}^n{n\choose k}{n+k\choose k}c_k^{(r)},
\end{equation*}
Is it true that all the numbers $c_k^{(r)}$ are integers?
\end{problem}

At the end of his paper, Zudilin~\cite{Zu} raised the problem
of finding and solving a $q$-analogue of Problem~\ref{prob:Sch}.
Our second objective is to provide such a $q$-analogue.

{}For any integer $n$, define the \emph{$q$-shifted factorial} $(a)_n$ by
$(a)_0=1$ and
$$
(a)_n=
\begin{cases}
    (1-a)(1-aq)\cdots (1-aq^{n-1}), & \text{$n=1,2,\ldots,$} \\
    (1-aq^{-1})(1-aq^{-2})\cdots (1-aq^{n})^{-1}, & \text{$n=-1,-2,\ldots.$}
\end{cases}
$$
We will also use the compact notations for $m\geq 1$:
\begin{align*}
 (a_1,\ldots,a_m)_n:=(a_1)_n\cdots (a_m)_n,
\qquad
(a_1,\ldots,a_m)_\infty:=\lim_{n\to\infty}(a_1,\ldots,a_m)_n.
\end{align*}
The $q$-binomial coefficients are defined as
$$
{n\brack k}:={n\brack k}_q
=\frac{(q)_n}{(q)_k (q)_{n-k}}.
$$
Since $\frac{1}{(q)_n}=0$ if $n<0$, we have ${n\brack k}=0$ if $k>n$ or $k<0$.

The following is our first generalization of Calkin's result.
\begin{thm}\label{thm:qidentity}
For $m\geq 3$ and all positive integers $n_1,\ldots,n_m$, there holds
\begin{align}
&\sum_{k=-n_1}^{n_1}(-1)^k
q^{(m-1)k^2+{k\choose 2}}
\prod_{i=1}^m {n_i+n_{i+1}\brack n_i+k}\nonumber\\
&\hspace{2cm}
={n_1+n_m\brack n_1}\sum_{\lambda}
\prod_{i=1}^{m-2}q^{\lambda_i^2}{\lambda_{i-1}\brack \lambda_i}
{n_{i+1}+n_{i+2}\brack n_{i+1}-\lambda_i},  \label{eq:multi}
\end{align}
where $n_{m+1}=\lambda_0=n_1$ and the sum is over all  sequences
$\lambda=(\lambda_1,\ldots,\lambda_{m-2})$ of nonnegative integers
such that $\lambda_0\geq \lambda_1\geq \cdots \geq \lambda_{m-2}$.
\end{thm}

Calkin~\cite{Calkin} has given a partial $q$-analogue
of \eqref{eq:Calkin} by considering the alternating sum
$\sum_{k=0}^n(-1)^kq^{jk}{n\brack k}^m$.
In this respect,  besides \eqref{eq:multi}, we
shall also prove the following divisibility result.
\begin{thm}\label{thm:qbino}
For all positive integers $n_1,\ldots,n_m$, $n_{m+1}=n_1$,  the alternating sum
\begin{align*}
S(n_1,\ldots,n_m;j,q):={n_1+n_m\brack n_1}^{-1}\sum_{k=-n_1}^{n_1}(-1)^k
q^{jk^2+{k\choose 2}}
\prod_{i=1}^m {n_i+n_{i+1}\brack n_i+k}
\end{align*}
 is a polynomial in $q$ with nonnegative integral coefficients for
$0\leq j\leq m-1$.
\end{thm}

We shall give two proofs of Theorem~\ref{thm:qidentity}: The first one
is based on a recurrence relation formula for $S(n_1,\ldots,n_m;j,q)$,
which also  leads to a proof of Theorem~\ref{thm:qbino}. The second one
follows directly from Andrews' basic hypergeometric identity between a
single sum and a multiple
sum~\cite[Theorem 4]{Andrews75}.
\begin{thm}[Andrews~\cite{{Andrews75}}]\label{thm:andrews}
 For every integer $m\geq 0$, the following identity holds:
\begin{multline}\label{andrews}
\sum_{k\geq 0}\frac{(a,q\sqrt{a},-q\sqrt{a},b_1,c_1,\dots,b_m,c_m,q^{-N})_k}
{(q,\sqrt{a},-\sqrt{a},aq/b_1,aq/c_1,\dots,aq/b_m,aq/c_m,aq^{N+1})_k}
\left(\frac{a^mq^{m+N}}{b_1c_1\cdots b_mc_m}\right)^k\\
=\frac{(aq,aq/b_mc_m)_N}{(aq/b_m,aq/c_m)_N}\sum_{l_1,\dots,l_{m-1}\geq0}
\frac{(aq/b_1c_1)_{l_1}\cdots(aq/b_{m-1}c_{m-1})_{l_{m-1}}}{(q)_{l_1}\cdots(q)_{l_{m-1}}}\\
\times\frac{(b_2,c_2)_{l_1}\dots(b_m,c_m)_{l_1+\dots+l_{m-1}}}{(aq/b_1,aq/c_1)_{l_1}
\dots(aq/b_{m-1},aq/c_{m-1})_{l_1+\dots+l_{m-1}}}\\
\times\frac{(q^{-N})_{l_1+\dots+l_{m-1}}}{(b_mc_mq^{-N}/a)_{l_1+\dots+l_{m-1}}}
\frac{(aq)^{l_{m-2}+\dots+(m-2)l_1}
q^{l_1+\dots+l_{m-1}}}{(b_2c_2)^{l_1}\cdots(b_{m-1}c_{m-1})^{l_1+\dots+l_{m-2}}}.
\end{multline}
\end{thm}
It is interesting to note that  Theorem~\ref{thm:andrews}
is  also a key
ingredient in Zudilin's approach to  Problem~\ref{prob:Sch}.
Therefore, using Theorem~\ref{thm:andrews} in its full generality we are able to
 formulate and prove a $q$-analogue of Problem~\ref{prob:Sch}.
\begin{thm}\label{thm:qProb}
For any integer $r\geq 1$, define rational fractions $c^{(r)}_k(q)$
of the variable $q$,
independent of $n$, by writing
\begin{align}\label{qzudilin}
&\sum_{k=0}^nq^{r{n-k\choose 2}+(1-r){n\choose 2}}
{n\brack k}^r{n+k\brack k}^r
=\sum_{k=0}^nq^{{n-k\choose 2}+(1-r){k\choose 2}}
{n\brack k}{n+k\brack k}c^{(r)}_k(q).
\end{align}
Then $c^{(r)}_n(q)\in\mathbb{N}[q]$.
\end{thm}

Since the $r=1$ case is trivial, we may suppose that $r\geq 2$ in what follows.
As ${n\brack k}{n+k\brack k}={2k\brack k}{n+k\brack n-k}$,
invoking the
$q$-Legendre transform, which is
 a special case of  Carlitz's $q$-Gould-Hsu inverse formula~\cite{Ca} (see also \cite{Kr}):
\begin{align*}
a_n=\sum_{k=0}^nq^{{n-k\choose 2}}{n+k\brack n-k}b_k
\Longleftrightarrow b_n=\sum_{k=0}^n(-1)^{n-k}
\frac{1-q^{2k+1}}{1-q^{n+k+1}}{2n\brack n-k}a_k,
\end{align*}
we derive immediately from \eqref{qzudilin} that
\begin{equation*} 
q^{(1-r){n\choose 2}}{2n\brack n} c^{(r)}_n(q)
=\sum_{j=0}^n{2j\brack j}^r t_{n,j}^{(r)}(q),
\end{equation*}
where
\begin{equation*}
t_{n,j}^{(r)}(q)=q^{r{j+1\choose 2}}\sum_{k=j}^n(-1)^{n-k}
\frac{1-q^{2k+1}}{1-q^{n+k+1}}{2n\brack n-k}
{k+j\brack k-j}^r q^{{k\choose 2}-rjk}.
\end{equation*}
Therefore, Theorem~\ref{thm:qProb} is a consequence of the following theorem, which
 is our $q$-analogue of Zudilin's result \eqref{eq:strehl}.
\begin{thm}\label{thm:qZud}
For any integer $r\geq 2$, we have
\begin{equation*} 
q^{(r-1){n\choose 2}}{2j\brack j}{2n\brack n}^{-1}
t_{n,j}^{(r)}(q)\in\mathbb{N}[q].
\end{equation*}
\end{thm}
As will be shown, Theorem~\ref{thm:qZud}
follows directly from Andrews' identity~\eqref{andrews}.

This paper is organized as follows.
We will prove Theorems \ref{thm:qidentity}
and \ref{thm:qbino} in the next section.
The proof of Theorem~\ref{thm:qZud} is given in Section~3.
Some interesting divisibility results are given in Section~4.
In the last section we will present four related conjectures.
\section{Proof of Theorems \ref{thm:qidentity} and \ref{thm:qbino}}
We will need two known identities in $q$-series.
One is the $q$-Pfaff-Saalsch\"utz identity
\cite[Appendix (II.12)]{GR} (see also \cite{GZ05,Zeilberger}):
\begin{align}
{n_1+n_2\brack n_1+k}{n_2+n_3\brack n_2+k}{n_3+n_1\brack n_3+k}
=\sum_{r=0}^{n_1-k}\frac{q^{k^2+2kr} (q)_{n_1+n_2+n_3-k-r}}
{(q)_r (q)_{r+2k}(q)_{n_1-k-r}(q)_{n_2-k-r}(q)_{n_3-k-r}}, \label{eq:qpfaff}
\end{align}
where $\frac{1}{(q)_n}=0$ if $n<0$, and the other is the $q$-Dixon identity:
\begin{align}
\sum_{k=-n_1}^{n_1}(-1)^k q^{(3k^2-k)/2}
{n_1+n_2\brack n_1+k}{n_2+n_3\brack n_2+k}{n_3+n_1\brack n_3+k}
=\frac{(q)_{n_1+n_2+n_3}}{(q)_{n_1}(q)_{n_2}(q)_{n_3}}. \label{eq:qdixon}
\end{align}
A short proof of \eqref{eq:qdixon} is given in \cite{GZdixon}.

We first establish the following
recurrence formula.
\begin{lem}\label{lem:s}
Let $m\geq 3$. Then for all positive integers $n_1,\ldots,n_m$ and any integer
$j$, the following
recurrence holds:
\begin{align}
S(n_1,\ldots,n_m;j,q)
=\sum_{l=0}^{n_1} q^{l^2}{n_1\brack l}{n_2+n_3\brack n_2-l}S(l,n_3,\ldots,n_m;j-1,q).
\label{eq:lem}
\end{align}
\end{lem}
\pf
For any integer $k$ and positive integers $a_1,\ldots,a_l$, let
$$
C(a_1,\ldots,a_l;k)=\prod_{i=1}^l {a_i+a_{i+1}\brack a_i+k},
$$
where $a_{l+1}=a_1$. Then
\begin{equation}\label{eq:rewriting}
S(n_1,\ldots,n_m;j,q)
=\frac{(q)_{n_1}(q)_{n_m}}{(q)_{n_1+n_m}}
\sum_{k=-n_1}^{n_1}(-1)^k q^{jk^2+{k\choose 2}}
C(n_1,\ldots,n_m;k).
\end{equation}

We observe  that  for $m\geq 3$, we have
$$
C(n_1,\ldots,n_m;k)=\frac{(q)_{n_2+n_3}(q)_{n_m+n_1}}{(q)_{n_1+n_2}(q)_{n_m+n_3}}
{n_1+n_2\brack n_1+k}{n_1+n_2\brack n_2+k}C(n_3,\ldots,n_m;k),
$$
and, by letting $n_3\to\infty$ in \eqref{eq:qpfaff},
$$
{n_1+n_2\brack n_1+k}{n_1+n_2\brack n_2+k}
=\sum_{r=0}^{n_1-k}\frac{q^{r^2+2kr}(q)_{n_1+n_2}}{(q)_r(q)_{r+2k}
(q)_{n_1-k-r}(q)_{n_2-k-r}}.
$$
Plugging these into \eqref{eq:rewriting} we can write its right-hand side as
$$
R:=\sum_{k=-n_1}^{n_1} \sum_{r=0}^{n_1-k}(-1)^k C(n_3,\ldots,n_m;k)
\frac{q^{(r+k)^2+(j-1)k^2+{k\choose 2}}(q)_{n_2+n_3}(q)_{n_1}(q)_{n_m}}
{(q)_r(q)_{r+2k}(q)_{n_1-k-r}(q)_{n_2-k-r}(q)_{n_m+n_3}}.
$$
Setting $l=r+k$, then $-n_1\leq l\leq n_1$, but if $l<0$, at least
one of the indices $l+k$ and $l-k$ is negative for any integer $k$,
which implies that  $\frac{1}{(q)_{l-k}(q)_{l+k}}=0$ by convention.
Therefore,  exchanging the order of summation, we have
$$
R=\sum_{l=0}^{n_1} \frac{q^{l^2}(q)_{n_2+n_3}(q)_{n_1}(q)_{n_m}}
{(q)_{n_1-l}(q)_{n_2-l}(q)_{n_m+n_3}}
\sum_{k=-l}^l (-1)^k C(n_3,\ldots,n_m;k)
\frac{q^{(j-1)k^2+{k\choose 2}}}{(q)_{l-k}(q)_{l+k}}.
$$
Now, in the last sum making the substitution
$$
C(n_3,\ldots, n_m;k)
=\frac{(q)_{l-k}(q)_{l+k}(q)_{n_m+n_3}}{(q)_{n_3+l}(q)_{n_m+l}}C(l,n_3,\ldots, n_m;k),
$$
we obtain
the right-hand side of \eqref{eq:lem}. \qed

\noindent \emph{First proof of Theorem~\ref{thm:qidentity}}.
Letting $n_3\to\infty$ in \eqref{eq:qdixon} yields that
\begin{align}
S(n_1,n_2;1,q)=1. \label{eq:p11mn}
\end{align}
Theorem~\ref{thm:qidentity} then follows
by iterating $(m-2)$ times
formula~\eqref{eq:lem}. \qed

\medskip

\noindent
\noindent\emph{Second proof of Theorem~\ref{thm:qidentity}}.
Since
$$
{M\brack N+k}=(-1)^kq^{(M-N)k-{k\choose 2}}{M\brack N}\frac{(q^{-M+N})_k}
{(q^{N+1})_k},
$$
by collecting the terms of index $k$ and $-k$,
the left-hand side of \eqref{eq:multi} can be written as
\begin{align*}
L&:=\prod_{i=1}^m{n_i+n_{i+1}\brack n_i}+\sum_{k=1}^{n_1}(1+q^k)(-1)^k
q^{(m-1)k^2+{k\choose 2}}
\prod_{i=1}^m{n_i+n_{i+1}\brack n_i+k}\\
&= \prod_{i=1}^m{n_i+n_{i+1}\brack n_i}\left\{1+
\sum_{k=1}^{n_1}(1+q^k)(-1)^{(m-1)k}
q^{(m-1){k+1\choose 2}}
\prod_{i=1}^mq^{n_ik}\frac{(q^{-n_{i+1}})_k}{(q^{n_{i}+1})_k}\right\}.
\end{align*}
Letting $c_1=c_2=\cdots =c_m=c\to \infty$ and $a\to 1$ in Andrews'
formula~\eqref{andrews}  we get
\begin{align}
1+&\sum_{k\geq 1}(1+q^k)\frac{(b_1,\ldots,b_m, q^{-N})_k}
{(q/b_1,\ldots,q/b_m, q^{N+1})_k}
(-1)^{mk}q^{m{k\choose 2}}
\left(\frac{q^{m+N}}{b_1b_2\ldots b_m}\right)^k\nonumber\\
&=\frac{(q)_N}{(q/b_m)_N}\sum_{l_1,\ldots, l_{m-1}\geq 0}
\frac{(q)_N}{(q)_{l_1}\cdots (q)_{l_{m-1}}(q)_{N-l_1-\cdots- l_{m-1}}}
\nonumber\\
&\times\prod_{i=1}^{m-1}\frac{(b_{i+1})_{l_1+\cdots+l_{i}}}
{(q/b_{i})_{l_1+\cdots +l_{i}}}
\left(\frac{-1}{b_{i+1}}\right)^{l_1+\cdots+l_{i}}
q^{{l_1+\cdots+l_{i}\choose 2}+(m-i)l_i}.\label{eq:andrewslimit}
\end{align}
Now, shifting $m$ to $m-1$ in \eqref{eq:andrewslimit}, setting
$$
N=n_m,\quad  b_i=q^{-n_i} \quad \textrm{for}
\quad i=1,\ldots, m-1,
$$
and $\lambda_i=l_1+\cdots +l_i$ for $i=1,\ldots, m-2$, one sees that
$L$ equals
\begin{align*}
&\hskip -2mm
\prod_{i=1}^m{n_i+n_{i+1}\brack n_i}
\frac{(q)_{n_m}^2}{(q^{1+n_{m-1}})_{n_m}}
\sum_{0\leq \lambda_1\leq\cdots\leq \lambda_{m-2}}
\prod_{i=1}^{m-2}\frac{(q^{-n_{i+1}})_{\lambda_{i}}
(-1)^{\lambda_{i}}q^{{\lambda_{i}+1\choose 2}+n_{i+1}\lambda_i}}
{(q^{1+n_i})_{\lambda_{i}}(q)_{\lambda_i-\lambda_{i-1}}}
\nonumber\\
&={n_1+n_m\brack n_1}\sum_{0\leq \lambda_1\leq\cdots\leq \lambda_{m-2}}
\prod_{i=1}^{m-2}q^{\lambda_i^2}{\lambda_{i+1}\brack \lambda_i}
{n_i+n_{i+1}\brack n_i+\lambda_i},
\end{align*}
where $\lambda_0=0$ and $\lambda_{m-1}=n_m$. The latter identity
 is clearly equivalent to Theorem~\ref{thm:qidentity}. \qed

In order to prove Theorem~\ref{thm:qbino}, we shall need the
following relation:
\begin{align}
S(n_1,\ldots,n_{m};0,q)
=S(n_1,\ldots,n_{m};m-1,q^{-1})\,q^{n_1n_2+n_2n_3+\cdots+n_{m-1}n_m}.
\label{eq:pinv}
\end{align}
As ${n\brack k}_{q^{-1}}={n\brack k} q^{k(k-n)}$,
Eq.~\eqref{eq:pinv} can be verified by substituting
$q$ by $q^{-1}$ and then replacing $k$ by $-k$ in the definition of
$S(n_1,\ldots,n_{m};m-1,q)$.

\medskip
\noindent
{\it Proof of Theorem~\ref{thm:qbino}.}
We proceed by induction on $m\geq 1$.
By the $q$-binomial theorem \cite[(II.3)]{GR},
we have
$$
S(n_1;0,q)=\sum_{k=-n_1}^{n_1}(-1)^kq^{k \choose 2}{2n\brack n+k}=0.
$$
In view
of \eqref{eq:p11mn}, it follows from \eqref{eq:pinv} that
\begin{align*}
S(n_1,n_2;0,q)=S(n_1,n_2;1,q^{-1})\,q^{n_1 n_2}=q^{n_1 n_2}.
\end{align*}
So the theorem is valid for $m\leq 2$.

Now suppose that
the expression $S(n_1,\ldots,n_{m-1};j,q)$
is a polynomial in $q$ with nonnegative integral coefficients for some $m\geq 3$ and $0\leq j\leq m-2$.
Then by the recurrence formula \eqref{eq:lem},
so is $S(n_1,\ldots,n_{m};j,q)$ for $1\leq j\leq m-1$.
It remains to show that $S(n_1,\ldots,n_{m};0,q)$ has the required
property. By Theorem~\ref{thm:qidentity} we know that $S(n_1,\ldots,n_{m};m-1,q)$
is a polynomial in $q$.
Since  the $q$-binomial coefficient ${n\brack k}$ is a polynomial in $q$ of degree
$k(n-k)$ (see \cite[p.~33]{Andrews98}), it is easy to see from the definition
of $S(n_1,\ldots,n_{m};m-1,q)$ that the degree of the polynomial
$S(n_1,\ldots,n_{m};m-1,q)$ is less than or equal to
$n_1n_2+n_2n_3+\cdots+n_{m-1}n_m$.
It follows from \eqref{eq:pinv} that $S(n_1,\ldots,n_{m};0,q)$
is also a polynomial in $q$ with nonnegative integral coefficients.
This completes the inductive step of the proof.
\qed

\begin{rmk}
Though it is not necessary to check the $m=3$ case to valid our induction argument, we think it is convenient
to include here
the formulas for $m=3$. First, the $q$-Dixon identity \eqref{eq:qdixon} implies that
\begin{align*}
S(n_1,n_2,n_3;1,q)={n_1+n_2+n_3\brack n_2}.
\end{align*}
{}From \eqref{eq:lem} and \eqref{eq:p11mn} we derive
\begin{align*}
S(n_1,n_2,n_3;2,q)=\sum_{l=0}^{n_1} q^{l^2}{n_1\brack l}{n_2+n_3\brack n_2-l}.
\end{align*}
Finally, applying \eqref{eq:pinv} we get
\begin{align*}
S(n_1,n_2,n_3;0,q)
&=S(n_1,n_2,n_3;2,q^{-1})\,q^{n_1n_2+n_2n_3}  \\
&=\sum_{l=0}^{n_1} q^{(n_1-l)(n_2-l)+n_3l}{n_1\brack l}{n_2+n_3\brack n_2-l}.
\end{align*}
\end{rmk}
\section{Proof of Theorem~\ref{thm:qZud} }
We will
distinguish the cases where $r\geq 2$ is even or odd,
and treat separately the values $r=2$ and $r=3$.
\begin{itemize}
\item For $r=2$, apply (\ref{andrews}) specialized with $m=1$, $a=q^{-(2n+1)}$,
$N=n-j$, $b_1=q^{-n}$ and $c_1=q^{-(n-j)}$.
The left-hand side of (\ref{andrews}) is then equal to
$$
{n+j\brack 2j}^{-2} q^{-2{n-j\choose 2}+{n\choose 2}} t_{n,j}^{(2)}(q).
$$
Equating this with the right-hand side gives
$$
t_{n,j}^{(2)}(q)=\frac{(q)_{2n}(q)_j^2}{(q)_n(q)_{2j}(q)_{2j-n}(q)_{n-j}^2}
q^{2{n-j\choose 2}-{n\choose 2}},
$$
which shows that ${2j\brack j}{2n\brack n}^{-1}
q^{{n\choose 2}} t_{n,j}^{(2)}(q)\in\mathbb{N}[q]$.\\

\item
For $r=3$, apply (\ref{andrews}) specialized with $m=1$, $a=q^{-(2n+1)}$, $N=n-j$
and $b_1=c_1=q^{-(n-j)}$. This yields in that case
$$
t_{n,j}^{(3)}(q)=\frac{(q)_{2n}}{(q)_{3j-n}(q)_{n-j}^3}
q^{3{n-j\choose 2}-2{n\choose 2}},
$$
which shows that ${2j\brack j}{2n\brack n}^{-1}
q^{2{n\choose 2}} t_{n,j}^{(3)}(q)\in\mathbb{N}[q]$.\\

\item For $r=2s\geq 4$, apply (\ref{andrews}) with $m=s\geq 2$,
$a=q^{-(2n+1)}$, $N=n-j$, $b_1=q^{-n}$ and $c_1=b_i=c_i=q^{-(n-j)}$,
$\forall\, i\in\{2,\dots,s\}$ to get
\begin{multline*}
q^{(2s-1){n\choose 2}-2s{n-j\choose 2}}\,t_{n,j}^{(2s)}(q)=\\
\frac{(q)_{2n}(q)_j}{(q)_n(q)_{2j}(q)_{n-j}}\sum_{l_1\geq0}{j\brack l_1}
{n-l_1\brack j}{n-l_1+j\brack n-l_1-j}
q^{{l_1\choose 2}+2j(s-1)l_1+(j+1-n)l_1}\\
\times\sum_{l_2\geq0}{2j\brack l_2}
{n-l_1-l_2+j\brack n-l_1-l_2-j}^2
q^{{l_2\choose 2}+2j(s-2)l_2+(j+1-n)l_2}\times\cdots\\
\times\sum_{l_{s-1}\geq0}{2j\brack l_{s-1}}
{n-l_1-\dots-l_{s-1}+j\brack n-l_1-\dots-l_{s-1}-j}^2
q^{{l_{s-1}\choose 2}+2jl_{s-1}+(j+1-n)l_{s-1}}\\
\times{2j\brack n-l_1-\dots-l_{s-1}-j}
q^{{l_1+\dots+l_{s-1}\choose 2}}.
\end{multline*}
As the condition $ l_1+\dots+l_{s-1}\leq n-j$ holds in the last summation, we can see that for $s\geq2$,
$
{2j\brack j}{2n\brack n}^{-1}
q^{(2s-1){n\choose 2}}\,t_{n,j}^{(2s)}(q)\in\mathbb{N}[q].
$

\item For $r=2s+1\geq 5$, apply (\ref{andrews}) with
$m=s\geq 2$, $a=q^{-(2n+1)}$, $N=n-j$, and $b_i=c_i=q^{-(n-j)}$,
$\forall\, i\in\{1,\dots,s\}$ to get
\begin{multline*}
q^{2s{n\choose 2}-(2s+1){n-j\choose 2}}t_{n,j}^{(2s+1)}(q)=\\
\frac{(q)_{2n}}{(q)_{2j}(q)_{n-j}^2}\sum_{l_1\geq0}{2j\brack l_1}
{n-l_1+j\brack n-l_1-j}^2q^{{l_1\choose 2}+2j(s-1)l_1+(j+1-n)l_1}\\
\times\sum_{l_2\geq0}{2j\brack l_2}
{n-l_1-l_2+j\brack n-l_1-l_2-j}^2
q^{{l_2\choose 2}+2j(s-2)l_2+(j+1-n)l_2}\times\cdots\\
\times\sum_{l_{s-1}\geq0}{2j\brack l_{s-1}}
{n-l_1-\dots-l_{s-1}+j\brack n-l_1-\dots-l_{s-1}-j}^2
q^{{l_{s-1}\choose 2}+2jl_{s-1}+(j+1-n)l_{s-1}}\\
\times{2j\brack n-l_1-\dots-l_{s-1}-j}q^{{l_1+\dots+l_{s-1}\choose 2}}.
\end{multline*}
As the condition $l_1+\dots+l_{s-1}\leq n-j$ holds in the last summation,
we can see that for $s\geq 2$,
$
{2j\brack j}{2n\brack n}^{-1}
q^{2s{n\choose 2}}\,t_{n,j}^{(2s+1)}(q)\in\mathbb{N}[q].
$
\end{itemize}
\begin{rmk}
In the special case $r=2$, our proof gives the following expression
for the coefficients $c^{(2)}_n(q)$:
\begin{equation}\label{r=2}
c_n^{(2)}(q)=\sum_{j=0}^n{2j\brack n}{n\brack
j}^2q^{2{n-j\choose 2}}.
\end{equation}
These coefficients are $q$-analogues of the famous $c^{(2)}_n(1)$ involved
in Apery's proof of the irrationality of $\zeta(3)$:
\begin{equation}\label{r=2q=1}
c_n^{(2)}(1)=\sum_{j=0}^n{2j\choose n}{n\choose j}^2=\sum_{j=0}^n{n\choose j}^3.
\end{equation}
As explained in \cite{St}, when $q=1$,
one can derive the last expression from (\ref{r=2}) in an
elementary way (by two iteration of the Chu-Vandermonde formula).
But our $q$-analogue (\ref{r=2}) does not
lead to a natural $q$-analogue of (\ref{r=2q=1}).
\end{rmk}

\section{Consequences of Theorems \ref{thm:qidentity} and \ref{thm:qbino}}
Letting $q=1$ in Theorem~\ref{thm:qidentity} we obtain a
direct generalization of Calkin's result~\eqref{eq:Calkin}.
\begin{thm}\label{thm:bino}
For $m\geq 3$ and all positive integers $n_1,\ldots,n_m$, there holds
\begin{equation}\label{eq:calkingeneral}
\sum_{k=-n_1}^{n_1}(-1)^k\prod_{i=1}^m
{n_i+n_{i+1}\choose n_i+k}
={n_1+n_m\choose n_1}\sum_{\lambda}
\prod_{i=1}^{m-2}{\lambda_{i-1}\choose \lambda_i}{n_{i+1}+n_{i+2}
\choose n_{i+1}-\lambda_i},
\end{equation}
 where $n_{m+1}=\lambda_0=n_1$ and the sum is over all  sequences
$\lambda=(\lambda_1,\ldots,\lambda_{m-2})$ of nonnegative integers
such that $\lambda_0\geq \lambda_1\geq \cdots \geq \lambda_{m-2}$.
\end{thm}

\begin{rmk}
For $m=1$ and $2$, it is easy to see that the left-hand side of
\eqref{eq:calkingeneral} is equal to 0 and ${n_1+n_2\choose n_1}$,
respectively. Calkin's result follows from \eqref{eq:calkingeneral}
by setting $n_i=n$ for all $i=1,\ldots, m$.
\end{rmk}

Letting $n_1=\cdots=n_m=n$ in Theorem~\ref{thm:qbino}, we obtain a complete
$q$-analogue of Calkin's result.
\begin{cor}
For all positive $m$, $n$ and $0\leq j\leq m-1$,
$$
{2n\brack n}^{-1}\sum_{k=-n}^n (-1)^k q^{jk^2+{k\choose 2}} {2n\brack n+k}^m
$$
is a polynomial in $q$ with nonnegative integral coefficients.
\end{cor}

Letting $n_{2i-1}=m$ and $n_{2i}=n$ for $1\leq i\leq r$ in Theorem~\ref{thm:qbino}, we obtain
\begin{cor}
For all positive $m$, $n$, $r$ and $0\leq j\leq 2r-1$,
$$
{m+n\brack m}^{-1}\sum_{k=-m}^m (-1)^k q^{jk^2+{k\choose 2}}
{m+n\brack m+k}^r {m+n\brack n+k}^r
$$
is a polynomial in $q$ with nonnegative integral coefficients. In particular,
$$
\sum_{k=-m}^m (-1)^k {m+n\choose m+k}^r {m+n\choose n+k}^r
$$
is divisible by ${m+n\choose m}$.
\end{cor}

Letting $n_{3i-2}=l$, $n_{3i-1}=m$ and $n_{3i}=n$ for $1\leq i\leq r$
in Theorem~\ref{thm:qbino}, we obtain
\begin{cor}
For all positive $l$, $m$, $n$, $r$ and $0\leq j\leq 3r-1$,
$$
{l+n\brack n}^{-1}\sum_{k=-l}^l (-1)^k q^{jk^2+{k\choose 2}}
{l+m\brack l+k}^r {m+n\brack m+k}^r {n+l\brack n+k}^r
$$
is a polynomial in $q$ with nonnegative integral coefficients. In particular,
$$
\sum_{k=-l}^l (-1)^k {l+m\choose l+k}^r{m+n\choose m+k}^r {n+l\choose n+k}^r
$$
is divisible by ${l+m\choose l}$, ${m+n\choose m}$ and ${n+l\choose n}$.
\end{cor}

Letting $m=2r+s$, $n_1=n_3=\cdots=n_{2r-1}=n+1$ and let all the other $n_i$ be $n$ in
Theorem~\ref{thm:bino}, we get
\begin{cor}\label{cor:nn+1}
For all positive $r$, $s$ and $n$,
\begin{align*}
\sum_{k=-n}^n (-1)^k {2n+1\choose n+k+1}^r{2n+1\choose n+k}^r{2n\choose n+k}^s
\end{align*}
is divisible by both ${2n\choose n}$ and ${2n+1\choose n}$, and is therefore
divisible by $(2n+1){2n\choose n}$.
\end{cor}

However, the following result is not a special case of Theorem~\ref{thm:bino}.
\begin{cor}
For all nonnegative $r$ and $s$ and positive $t$ and $n$,
\begin{align*}
\sum_{k=-n}^n
(-1)^k {2n+1\choose n+k+1}^r{2n+1\choose n+k}^s{2n\choose n+k}^t.
\end{align*}
is divisible by ${2n\choose n}$.
\end{cor}

\pf We proceed by induction on $|r-s|$. The $r=s$ case is clear from
Corollary~\ref{cor:nn+1}. Suppose the statement is true for $|r-s|\leq m-1$.
By Theorem~\ref{thm:bino}, one sees that
\begin{align}
&\hskip -3mm
\sum_{k=-n}^n(-1)^k{2n+2\choose n+k+1}^{m}{2n+1\choose n+k+1}^s{2n+1\choose n+k}^s
{2n\choose n+k}^t \nonumber\\
&=\frac{2n+2}{2n+1}\sum_{k=-n}^n(-1)^k{2n+2\choose n+k+1}^{m-1}
{2n+1\choose n+k+1}^{s+1}{2n+1\choose n+k}^{s+1}
{2n\choose n+k}^{t-1}, \label{eq:rsst}
\end{align}
where $m,t\geq 1$, is divisible by
$$
\frac{2n+2}{2n+1}{2n+1\choose n}=2{2n\choose n}.
$$
By the binomial theorem, we have
\begin{align}
{2n+2\choose n+k+1}^{m}
&=\left({2n+1\choose n+k+1}+{2n+1\choose n+k}\right)^m \nonumber\\
&=\sum_{i=0}^m{m\choose i}{2n+1\choose n+k+1}^i{2n+1\choose n+k}^{m-i}.
\label{eq:2n+2}
\end{align}
Substituting \eqref{eq:2n+2} into the left-hand of \eqref{eq:rsst}
and using the induction hypothesis and symmetry, we find that
\begin{align*}
&\hskip -3mm
\sum_{k=-n}^n (-1)^k {2n+1\choose n+k+1}^{m+s}{2n+1\choose n+k}^s{2n\choose n+k}^t  \\
&{}+\sum_{k=-n}^n (-1)^k {2n+1\choose n+k+1}^{s}{2n+1\choose n+k}^{m+s}{2n\choose n+k}^t
\end{align*}
is divisible by $2{2n\choose n}$. However, replacing $k$ by $-k$, one sees that
\begin{align*}
&\hskip -3mm
\sum_{k=-n}^n (-1)^k {2n+1\choose n+k+1}^{m+s}{2n+1\choose n+k}^s{2n\choose n+k}^t  \\
&=\sum_{k=-n}^n (-1)^k {2n+1\choose n+k+1}^{s}{2n+1\choose n+k}^{m+s}{2n\choose n+k}^t.
\end{align*}
This proves that the statement is true for $|r-s|=m$.  \qed

It is clear that Theorems \ref{thm:qbino} and \ref{thm:bino} can be
restated in the following forms.

\begin{thm}
For all positive integers $n_1,\ldots,n_m$ and $0\leq j\leq m-1$, the alternating sum
$$
(q)_{n_1}\prod_{i=1}^m\frac{(q)_{n_i+n_{i+1}}}{(q)_{2n_i}}
\sum_{k=-n_1}^{n_1}(-1)^k q^{jk^2+{k\choose 2}}\prod_{i=1}^m {2n_i\brack n_i+k},
$$
where $n_{m+1}=0$, is a polynomial in $q$ with nonnegative integral coefficients.
\end{thm}

\begin{thm}\label{thm:rebino}
For all positive integers $n_1,\ldots,n_m$, we have
$$
n_1!\prod_{i=1}^m\frac{(n_i+n_{i+1})!}{(2n_i)!}
\sum_{k=-n_1}^{n_1}(-1)^k\prod_{i=1}^m {2n_i\choose n_i+k}\in \N,
$$
where $n_{m+1}=0$.
\end{thm}

It is easy to see that, for all positive integers $m$ and $n$, the expression
$\frac{(2m)!(2n)!}{(m+n)!m!n!}$ is an integer. Letting $n_1=\cdots n_r=m$ and
$n_{r+1}=\cdots=n_{r+s}=n$ in Theorem~\ref{thm:rebino}, we obtain
\begin{cor}\label{cor:mn}
For all positive $m$, $n$, $r$ and $s$,
$$
\sum_{k=-m}^m(-1)^k{2m\choose m+k}^r{2n\choose n+k}^s
$$
is divisible by $\frac{(2m)!(2n)!}{(m+n)!m!n!}$.
\end{cor}

In particular, we find that
$$
\sum_{k=-n}^n (-1)^k {4n\choose 2n+k}^r{2n\choose n+k}^s
$$
is divisible by ${4n\choose n}$, and
$$
\sum_{k=-n}^n (-1)^k {6n\choose 3n+k}^r{2n\choose n+k}^s
$$
is divisible by
$\frac{(6n)!(2n)!}{(4n)!(3n)!n!}$.

{}From Theorem~\ref{thm:rebino} it is easy to see that
$$
n_1!\prod_{i=1}^m\frac{(n_i+n_{i+1})!}{(2n_i)!}
\sum_{k=-n_1}^{n_1}(-1)^k\prod_{i=1}^m {2n_i\choose n_i+k}^{r_i},
$$
where $n_{m+1}=0$, is a nonnegative integer for all $r_1,\ldots,r_m\geq 1$.
For $m=3$, letting $(n_1,n_2,n_3)$ be $(n,3n,2n)$, $(2n,n,3n)$, or
$(2n,n,4n)$, we obtain the following two corollaries.
\begin{cor}\label{cor:rst-246n}
For all positive $r$, $s$, $t$ and $n$,
$$
\sum_{k=-n}^n (-1)^k {6n\choose 3n+k}^r {4n\choose 2n+k}^s
{2n\choose n+k}^t
$$
is divisible by both ${6n\choose n}$ and ${6n\choose 3n}$.
\end{cor}

\begin{cor}\label{cor:rst-248n}
For all positive $r$, $s$, $t$ and $n$,
$$
\sum_{k=-n}^n (-1)^k {8n\choose 4n+k}^r {4n\choose 2n+k}^s
{2n\choose n+k}^t
$$
is divisible by ${8n\choose 3n}$.
\end{cor}

\section{Some open problems}
Based on computer experiments, we would like to present four
interesting conjectures. The first two are refinements of
Corollaries~\ref{cor:rst-246n} and \ref{cor:rst-248n} respectively.
\begin{conj}\label{conj:5.1}
For all positive $r$, $s$, $t$ and $n$,
$$
\sum_{k=-n}^n (-1)^k {6n\choose 3n+k}^r {4n\choose 2n+k}^s
{2n\choose n+k}^t
$$
is divisible by both $2{6n\choose n}$ and $6{6n\choose 3n}$.
\end{conj}

\begin{conj} \label{conj:5.2}
For all positive $r$, $s$, $t$ and $n$ with $(r,s,t)\neq (1,1,1)$,
$$
\sum_{k=-n}^n (-1)^k {8n\choose 4n+k}^r {4n\choose 2n+k}^s
{2n\choose n+k}^t
$$
is divisible by $2{8n\choose 3n}$.
\end{conj}

Conjectures~\ref{conj:5.1} and \ref{conj:5.2} are true for
$r+s+t\leq 10$ and $n\leq 100$.

Let $\gcd(a_1,a_2,\ldots)$ denote the greatest common divisor of integers
$a_1,a_2,\ldots.$
\begin{conj}
For all positive $m$ and $n$, we have
\begin{align}
\gcd\left(\sum_{k=-n}^n (-1)^k{2n\choose n+k}^r,\ r=m,m+1,\ldots\right)
={2n\choose n}.
\end{align}
\end{conj}

Let $d\geq 2$ be a fixed integer. Every nonnegative integer $n$ can be uniquely written as
$$
n=\sum_{i\geq 0} a_id^i,
$$
where $0\leq a_i\leq d-1$ for all $i$ and only finitely many number of $b_i$ are nonzeros,
denoted by $n=[\cdots a_1a_0]_d$, in which the first $0$'s are omitted.
Let $n=[a_1\cdots a_r]_3=[b_1\cdots b_s]_7=[c_1\cdots c_t]_{13}$.
We now define three statistics $\alpha(n)$, $\beta(n)$ and $\gamma(n)$ as follows.
\begin{itemize}
\item Let $\alpha(n)$ be the number of disconnected $2$'s in the sequence
$a_1\cdots a_r$. Here two nonzero digits $a_i$ and $a_j$ are said to be
disconnected if there is at least one $0$ between $a_i$ and $a_j$.
\item Let $\beta(n)$ be the number of $1$'s in $b_1\cdots b_s$
which are not immediately followed by a $4$, $5$, or $6$.
\item Let $\gamma(n)$ be the number of $1$'s in $c_1\cdots c_t$ which are immediately
 followed by one of $7,\ldots,12$, or immediately followed by a number of $6$'s
and then followed by one of $7,\ldots,12$.
\end{itemize}
For instance, $[20212]_3=185$ and so $\alpha(185)=2$; $[10142]_7=2480$, and so $\beta(2480)=1$;
$[1667]_{13}=3296$ and so $\gamma(3296)=1$.
The first $n$ such that $\alpha(n)=4$ is $[2020202]_3=1640$; the first $n$ such that
$\beta(n)=4$ is $[1111]_7=400$;
while the first $n$ such that $\gamma(n)=4$ is $[17171717]_{13}=97110800$.

We end this paper with the following conjecture.
\begin{conj}
For every positive integer $n$, we have
\begin{align*}
\gcd\left(\sum_{k=-n}^n (-1)^k{2n\choose n+k}^{3r},\ r=1,2,\ldots\right)
&={2n\choose n}3^{\alpha(n)}, \\
\gcd\left(\sum_{k=-n}^n (-1)^k{2n\choose n+k}^{3r+1},\ r=1,2,\ldots\right)
&={2n\choose n}7^{\beta(n)}13^{\gamma(n)}, \\
\gcd\left(\sum_{k=-n}^n (-1)^k{2n\choose n+k}^{3r+2},\ r=1,2,\ldots\right)
&={2n\choose n}.
\end{align*}
\end{conj}

\section*{Acknowledgments} This work was partially done during the first
author's visit to Institut Camille Jordan, Universit\'e Claude
Bernard (Lyon I), and was supported by a French postdoctoral
fellowship. We thank Christian Krattenthaler for conveying us his
feeling that Theorem~\ref{thm:qidentity} should follow from Andrews'
formula~\cite[Theorem~4]{Andrews75}.
We also thank Wadim Zudilin for useful conversations during his visit in Lyon.
The first author is grateful to Hong-Xing Ding for helpful comments on
Conjecture~\ref{conj:5.2}.

\renewcommand{\baselinestretch}{1}

\noindent{Victor J. W. Guo$^1$, Fr\'ed\'eric Jouhet$^2$ and Jiang Zeng$^3$}
\vskip 2mm

{\footnotesize\noindent
$^1$Department of Mathematics, East China Normal University,
Shanghai 200062,\\
People's Republic of China\\
{\tt jwguo@math.ecnu.edu.cn,\quad http://math.ecnu.edu.cn/\textasciitilde{jwguo}}\\

\noindent $^2$Universit\'e de Lyon,
Universit\'e Lyon1,
CNRS, UMR 5208 Institut Camille Jordan,\\
B\^atiment du Doyen Jean Braconnier,
43, blvd du 11 novembre 1918,
F-69622 Villeurbanne Cedex, France\\
{\tt jouhet@math.univ-lyon1.fr,\quad http://math.univ-lyon1.fr/\textasciitilde{jouhet}}\\

\noindent $^3$Universit\'e de Lyon,
Universit\'e Lyon1,
CNRS, UMR 5208 Institut Camille Jordan,\\
B\^atiment du Doyen Jean Braconnier,
43, blvd du 11 novembre 1918,
F-69622 Villeurbanne Cedex, France\\
{\tt zeng@math.univ-lyon1.fr,\quad http://math.univ-lyon1.fr/\textasciitilde{zeng}}
}

\end{document}